\documentclass[12pt,reqno]{amsart}
\usepackage{amssymb}
\usepackage{mathrsfs}
\newtheorem{theorem}{Theorem}
\newtheorem*{theorem*}{Theorem}
\newtheorem{lemma}{Lemma}
\newtheorem{prop}{Proposition}
\newtheorem{corollary}{Corollary}
\theoremstyle{definition}
\newtheorem{definition}{Definition}
\newtheorem*{definition*}{Definition}
\theoremstyle{remark}
\newtheorem*{example*}{Example}
\newtheorem{remark}{Remark}
\newcommand*{\ptens}[1]{\mathop{\widehat\otimes}_{#1}}
\newcommand*{\Ptens}{\mathop{\widehat\otimes}}
\newcommand*{\lmod}{\mbox{-}\!\mathop{\mathsf{mod}}}
\newcommand*{\rmod}{\mathop{\mathsf{mod}}\!\mbox{-}}
\newcommand*{\CC}{\mathbb C}
\newcommand*{\N}{\mathbb N}
\newcommand*{\eps}{\varepsilon}
\begin{document}
\title[Strictly Flat Cyclic Fr\'echet Modules]%
{Strictly Flat Cyclic Fr\'echet Modules\\ and Approximate Identities}
\author{A. Yu. Pirkovskii}
\address{Department of Nonlinear Analysis and Optimization\\
Faculty of Science\\
Peoples' Friendship University of Russia\\
Mikluho-Maklaya 6\\
117198 Moscow\\
Russia}
\email{pirkosha@sci.pfu.edu.ru, pirkosha@online.ru}
\thanks{Partially supported by the RFBR grants 05-01-00982 and 05-01-00001,
by the President of Russia grant MK-2049.2004.1, and by the Royal Society/NATO
Postdoctoral Fellowship, grant RG.MATH 447256.}
\begin{abstract}
Let $A$ be a locally $m$-convex
Fr\'echet algebra. We give a necessary and sufficient condition
for a cyclic Fr\'echet $A$-module $X=A_+/I$ to be strictly flat,
generalizing thereby a criterion of Helemskii and Sheinberg \cite{X1}.
To this end, we introduce a notion of locally bounded approximate identity (a.i.),
and we show that $X$ is strictly flat if and only if the ideal $I$
has a right locally bounded a.i. An example is given of a
commutative locally $m$-convex Fr\'echet algebra that has
a locally bounded a.i., but does not have a bounded a.i.
On the other hand, we show that a quasinormable
locally $m$-convex Fr\'echet algebra has a locally bounded a.i.
if and only if it has a bounded a.i. Some applications
to amenable Fr\'echet algebras are also given.
\end{abstract}
\maketitle

\section{Introduction}

Flat Banach modules over Banach algebras were introduced
by Helemskii \cite{Hel_period} by analogy with pure algebra.
An important fact which explains the r\^ole of flat modules
in Functional Analysis is that flatness is closely related
to amenability. More exactly, amenable Banach algebras
can be characterized in terms of flat Banach modules
(see \cite[VII.2]{X1}; cf. also Section~\ref{sec:amen} below).

The definition of flat Banach module readily extends to Fr\'echet
modules over Fr\'echet algebras. By a {\em Fr\'echet algebra}
we mean a complete, Hausdorff, metrizable locally convex $\CC$-algebra.
If $A$ is a Fr\'echet algebra, then a {\em left Fr\'echet $A$-module} is a
Fr\'echet space $X$ together with the structure of a left $A$-module such
that the product map $A\times X\to X$ is continuous. Right Fr\'echet
$A$-modules are defined similarly. The category of all left
(respectively, right) Fr\'echet $A$-modules will be denoted by $A\lmod$
(respectively, $\rmod A$). Recall from \cite{X1} that a chain complex $Y_\bullet$
of Fr\'echet $A$-modules is {\em admissible} if it splits in the category
of Fr\'echet spaces, i.e., if it has a contracting homotopy consisting
of continuous linear maps.

\begin{definition*}[Helemskii]
A left Fr\'echet $A$-module $X$ is said to be {\em flat}
(respectively, {\em strictly flat}) if for each admissible complex
(respectively, for each exact complex) $Y_\bullet$ in $\rmod A$
the sequence $Y_\bullet\ptens{A}X$ is exact.
\end{definition*}

Suppose now that $X$ is a cyclic left Fr\'echet $A$-module,
i.e., $X\cong A_+/I$ for a closed left ideal $I\subset A_+$.
Is is natural to ask when $X$ is flat or strictly flat. This question
is important, for example, because of its connection with
amenability (see \cite[VII.2]{X1} and Section~\ref{sec:amen} below).

In the case of Banach algebras, the answer
is as follows (see \cite[VII.1.4]{X1}).

\begin{theorem*}[Helemskii, Sheinberg]
\label{thm:Hel}
Let $A$ be a Banach algebra, and let $I\subset A_+$ be a closed left
ideal. Then the following conditions are equivalent:

\begin{itemize}
\item[{\upshape (i)}] $A_+/I$ is strictly flat;
\item[{\upshape (ii)}] $I$ has a right b.a.i.
\end{itemize}

If, in addition, $I$ is weakly complemented in $A_+$
(i.e., if the annihilator of $I$ is complemented in $A_+^*$),
then {\upshape (i)} and {\upshape (ii)} are equivalent to

\begin{itemize}
\item[{\upshape (iii)}] $A_+/I$ is flat.
\end{itemize}
\end{theorem*}

\begin{remark}
In the original form of the Helemskii-Sheinberg theorem, the phrase
``$X$ is flat'' (respectively, strictly flat) actually means
``$X$ is flat (respectively, strictly flat) as a Banach $A$-module'',
i.e., for each admissible (respectively, exact) complex
$Y_\bullet$ of right {\em Banach} $A$-modules
the sequence $Y_\bullet\ptens{A} X$ is exact. However,
this automatically implies the exactness of $Y_\bullet\ptens{A} X$ for each
admissible (respectively, exact) complex $Y_\bullet$ of right {\em Fr\'echet} $A$-modules,
as the following proposition suggests.
\end{remark}

\begin{prop}
Let $A$ be a Banach algebra, and let $X$ be a left Banach $A$-module.
Suppose that $X$ is flat (respectively, strictly flat) when considered
as a Banach $A$-module. Then $X$ is flat (respectively, strictly flat) as
a Fr\'echet $A$-module.
\end{prop}

The aim of this paper is to generalize the Helemskii-Sheinberg theorem to
Fr\'echet-Arens-Michael algebras (i.e., to locally $m$-convex Fr\'echet algebras).

\begin{remark}
Related results were independently obtained by C.~Podara.
\end{remark}

\section{The main result}
Let $A$ be an Arens-Michael algebra, i.e., a complete locally
$m$-convex algebra. Recall (see \cite{Michael} or \cite{X2}) that
$A$ is isomorphic to a projective limit of Banach algebras. More exactly,
let $\{\|\cdot\|_\lambda : \lambda\in\Lambda\}$ be a directed family
of submultiplicative seminorms generating the topology of $A$.
Given $\lambda,\mu\in\Lambda$, we write $\lambda\prec\mu$ if $\| a\|_\lambda\le\| a\|_\mu$
for all $a\in A$.
For each $\lambda\in\Lambda$ we set $N_\lambda=\{ a\in A : \| a\|_\lambda=0\}$.
Since $\|\cdot\|_\lambda$ is submultiplicative,
we see that $N_\lambda$ is a two-sided ideal of $A$, so that $A/N_\lambda$
is an algebra in a natural way. Moreover, the seminorm $\|\cdot\|_\lambda$
determines a submultiplicative norm on $A/N_\lambda$. Hence
the completion of $A/N_\lambda$ with respect to this norm is a Banach algebra.
This algebra is denoted by $A_\lambda$ and is called the {\em accompanying}
Banach algebra of $A$ corresponding to the seminorm $\|\cdot\|_\lambda$.
The canonical homomorphism $A\to A_\lambda$, $a\mapsto a+N_\lambda$, will be denoted
by $\tau_\lambda$. If $\lambda,\mu\in\Lambda$ and $\lambda\prec\mu$, then
there is a unique continuous homomorphism $\tau^\mu_\lambda\colon A_\mu\to A_\lambda$
such that $\tau_\lambda=\tau^\mu_\lambda \tau_\mu$. The family
$\{\tau_\lambda : \lambda\in\Lambda\}$ determines a continuous
homomorphism from $A$ to the projective limit $\varprojlim\{ A_\lambda,\tau^\mu_\lambda\}$,
and the {\em Arens-Michael decomposition theorem} states that this homomorphism
is a topological algebra isomorphism. The situation described above is
usually expressed by the phrase {\em ``Let $A$ be an Arens-Michael
algebra, and let $A=\varprojlim A_\lambda$ be the Arens-Michael decomposition
of $A$''}.

For each $\lambda\in\Lambda$ the homomorphism $\tau_\lambda\colon A\to A_\lambda$
uniquely extends to a unital homomorphism $\tau_\lambda^+\colon A_+\to (A_\lambda)_+$,
and a similar statement is true for all the connecting homomorphisms
$\tau^\mu_\lambda\; (\lambda\prec\mu)$. Obviously, $A_+\cong\varprojlim (A_\lambda)_+$.
Given a closed left ideal $I\subset A_+$, we set $I_\lambda=\overline{\tau_\lambda^+(I)}\subset (A_\lambda)_+$.
Evidently, $I_\lambda$ is a closed left ideal of $(A_\lambda)_+$.
Note that $\tau^\mu_\lambda(I_\mu)\subset I_\lambda$ for each $\lambda\prec\mu$,
and that $I\cong\varprojlim I_\lambda$.

From now on, we suppose that $A$ is a Fr\'echet-Arens-Michael algebra.
Without loss of generality, we may assume that $(\Lambda,\prec)=(\N,\le)$.

\begin{theorem}
\label{thm:I_n}
The following conditions are equivalent:
\begin{itemize}
\item[{\upshape (i)}]
$A_+/I$ is strictly flat;
\item[{\upshape (ii)}]
for each $n\in\N$, $I_n$ has a right b.a.i.
\end{itemize}
\end{theorem}

\begin{remark}
The implication (i)$\Longrightarrow$(ii) of the above theorem
holds for any Arens-Michael algebra. However, we do not know whether
(ii) implies (i) without the metrizability condition.
\end{remark}

It is easy to see that if $I$ has a right b.a.i., then so does $I_n$
for each $n\in\N$. A natural question is then whether the converse is also
true. As we shall see later, this need not be the case in general, but this
is the case under some additional linear topological assumptions on $I$.
Another natural question is whether it is possible to formulate
condition (ii) intrinsically, i.e., without referring to the accompanying
Banach algebras. Let us start by answering the latter question.

\section{Locally bounded approximate identities}
The following theorem is well known in the case of Banach algebras
(see, e.g., \cite{Dales,DW,Palmer}), but it readily extends
to arbitrary topological algebras.

\begin{theorem*}
Let $A$ be a topological algebra. Then
\begin{itemize}
\item[{\upshape (i)}]
$A$ has a right a.i. if and only if for each finite subset $F\subset A$
and each $0$-neighbourhood $U\subset A$ there exists $b\in A$ such that
$a-ab\in U$ for all $a\in F$;
\item[{\upshape (ii)}]
$A$ has a right b.a.i. if and only if there exists a bounded subset $S\subset A$
such that for each finite subset $F\subset A$
and each $0$-neighbourhood $U\subset A$ there exists $b\in S$ such that
$a-ab\in U$ for all $a\in F$.
\end{itemize}
\end{theorem*}

For our purposes, it is convenient to reformulate the above theorem
in the language of seminorms.

\begin{theorem*}
Let $A$ be a locally convex topological algebra, and let
$\{\|\cdot\|_\lambda : \lambda\in\Lambda\}$ be a directed family
of seminorms generating the topology of $A$. Then
\begin{enumerate}
\item[{\upshape (i)}]
$A$ has a right a.i. if and only if for each finite subset $F\subset A$,
each $\lambda\in\Lambda$, and each $\eps>0$ there exists $b\in A$ such that
$\| a-ab\|_\lambda<\eps$ for all $a\in F$;
\item[{\upshape (ii)}]
$A$ has a right b.a.i. if and only if there exists a family $\{ C_\lambda : \lambda\in\Lambda\}$
of positive reals such that for each finite subset $F\subset A$,
each $\lambda\in\Lambda$, and each $\eps>0$ there exists $b\in A$ such that
\begin{enumerate}
\item[{\upshape (ii$_1$)}]
$\| a-ab\|_\lambda<\eps$ for all $a\in F$, and
\item[{\upshape (ii$_2$)}]
$\| b\|_\mu\le C_\mu$ for all $\mu\in\Lambda$.
\end{enumerate}
\end{enumerate}
\end{theorem*}

Now let us relax condition (ii$_2$) as follows.

\begin{definition}
Let $A$ be a locally convex topological algebra, and let
$\{\|\cdot\|_\lambda : \lambda\in\Lambda\}$ be a directed family
of seminorms generating the topology of $A$. We say that
{\em $A$ has a right locally bounded a.i.} if
there exists a family $\{ C_\lambda : \lambda\in\Lambda\}$
of positive reals such that for each finite subset $F\subset A$,
each $\lambda\in\Lambda$, and each $\eps>0$ there exists $b\in A$ such that
\begin{enumerate}
\item[{\upshape (1)}]
$\| a-ab\|_\lambda<\eps$ for all $a\in F$, and
\item[{\upshape (2)}]
$\| b\|_\lambda\le C_\lambda$.
\end{enumerate}
\end{definition}

\begin{remark}
It is easy to see that the above definition does not depend on the
choice of the defining family of seminorms.
\end{remark}

\begin{prop}
\label{prop:lbai}
Let $A$ be an Arens-Michael algebra, and let
$A=\varprojlim A_\lambda$ be the Arens-Michael decomposition of $A$.
Then the following conditions are equivalent:
\begin{itemize}
\item[{\upshape (i)}]
$A$ has a right locally bounded a.i.;
\item[{\upshape (ii)}]
for each $\lambda\in\Lambda$, $A_\lambda$ has a right b.a.i.
\end{itemize}
\end{prop}

Now we can reformulate Theorem~\ref{thm:I_n} in a more elegant way.

\begin{theorem}
\label{thm:strflat2}
Let $A$ be a Fr\'echet-Arens-Michael algebra, and let $I\subset A_+$ be a closed
left ideal. Then the following conditions are equivalent:
\begin{itemize}
\item[{\upshape (i)}]
$A_+/I$ is strictly flat;
\item[{\upshape (ii)}]
$I$ has a right locally bounded a.i.
\end{itemize}
\end{theorem}

\section{Quasinormable Fr\'echet algebras}
It is clear that if a locally convex algebra $A$ has a locally
bounded a.i., then it has a bounded a.i. It is natural to ask
whether the converse is true. Before answering this question, let us recall
the following definition.

\begin{definition*}[Grothendieck \cite{Groth_F_DF}]
A locally convex space $E$ is {\em quasinormable} if
for each $0$-neighbourhood $U\subset E$ there exists a $0$-neighbourhood
$V\subset E$ such that for each $\eps>0$ there exists a bounded set
$B\subset E$ such that $V\subset B+\eps U$.
\end{definition*}

Many naturally arising Fr\'echet spaces are
quasinornable. Clearly, all Banach spaces and all Schwartz spaces \cite{Groth_F_DF}
are quasinormable. It is also true that each quojection (i.e., the projective limit
of a sequence of Banach spaces and surjective mappings) is quasinormable \cite{Dier_Zarn}.
This implies that the space $C(X)$ (where $X$ is a locally compact
Hausdorff topological space, countable at infinity)
and, more generally, each Fr\'echet locally $C^*$-algebra is quasinormable.
Standard examples of non-quasinormable Fr\'echet spaces belong to the class
of K\"othe sequence spaces \cite{Groth_F_DF} (see also \cite[Section 27]{MV}).
Within the class of function spaces, a number of concrete examples were found
in \cite{BonTask}.

\begin{theorem}
Let $A$ be a quasinormable Fr\'echet-Arens-Michael algebra with
a right locally bounded a.i. Then $A$ has a right b.a.i.
\end{theorem}

\section{A counterexample}
Now we present an example of a commutative
Fr\'echet-Arens-Michael algebra with a locally bounded a.i., but
without a b.a.i. Together with Theorem~\ref{thm:strflat2}, this will show
that Helemskii's theorem does not extend {\em verbatim} to Fr\'echet-Arens-Michael
algebras.

A ``building block'' for our example is the following Banach algebra.
In the sequel, for each $i\in\N$ we set $e_i=(0,\ldots,0,1,0,\ldots)$,
where the single nonzero entry is in the $i$th slot.
It is easy to see that there is a unique continuous multiplication
on $\ell^1$ such that $e_i e_j=e_{\min\{ i,j\}}$ for all $i,j\in\N$.
The resulting Banach algebra will be denoted by $A_1$. Let us remark
that $A_1$ is topologically isomorphic to the sequence algebra $bv_0$
consisting of all sequences converging to $0$ and having bounded variation
(see \cite{Dales}). An explicit isomorphism $A_1\to bv_0$ is given
by $e_n\mapsto e_1+\ldots +e_n$ ($n\in\N$). This implies that
$\{ e_n : n\in\N\}$ is a b.a.i. for $A_1$.

We shall need the following generalization of $A_1$. Let $P$ be
a family of real-valued sequences such that $p_i\ge 1$ for all
$p\in P$ and all $i\in\N$. Suppose also that $P$ is directed, i.e.,
for each $p,q\in P$ there exists $r\in P$ such that $r_i\ge\max\{ p_i,q_i\}$
for all $i\in\N$. Then it is easy to see that there exists a unique
multiplication on the K\"othe space
$$
\lambda(P)=\Bigl\{ a=(a_i)_{i\in\N}\in\CC^\N :
\| a\|_p=\sum_i |a_i| p_i < \infty\;\forall p\in P\Bigr\}
$$
such that $e_i e_j=e_{\min\{ i,j\}}$ for all $i,j\in\N$.
Moreover, we have $\| ab\|_p\le \| a\|_p \| b\|_p$ for all
$a,b\in\lambda(P)$ and all $p\in P$, so that $\lambda(P)$ becomes
as Arens-Michael algebra with respect to the above multiplication.
Let us denote this algebra by $A(P)$.

Given $p\in P$, the accompanying Banach algebra $A(P)_p$ can be
described in much the same way as $A_1$ (see above), by replacing
$\ell^1$ with the weighted space
$\ell^1(p)=\{ x=(x_i) : \| x\|=\sum_i |x_i|p_i<\infty\}$.
Together with Proposition~\ref{prop:lbai}, this easily implies the
following.

\begin{lemma}
\label{lemma:suff}
{\upshape (i)}
Suppose that each sequence $p\in P$ has a bounded subsequence. Then
$A(P)$ has a locally bounded a.i.

{\upshape (ii)} Suppose that there exists an infinite increasing
sequence $\{ n_k : k\in\N\}$ of positive integers such that
the sequence $\{ p_{n_k} : k\in\N\}$ is bounded for each $p\in P$.
Then $A(P)$ has a b.a.i.
\end{lemma}

Unfortunately, we do not know whether condition (ii) is necessary
for $A$ to have a b.a.i. In order to formulate a necessary condition,
let us introduce some notation. Given $a\in A(P)$, we set $w(a)=\sum_i |a_i|$
(note that this number is finite because $p_i\ge 1$ for all $p\in P$
and all $i\in\N$).
If $a\ne 0$, then we set $\ell(a)=\min\{ k\in\N : a_k\ne 0\}$.

\begin{lemma}
\label{lemma:necess}
Suppose that $A(P)$ has a b.a.i. Then there exists a bounded sequence
$\{ x_n\}\subset A(P)$ such that
\begin{equation}
\label{bai}
\ell(x_n)\to\infty\quad\text{as}\quad n\to\infty,\quad
\text{and}\quad \inf_n w(x_n)>0.
\end{equation}
\end{lemma}

Thus our aim is to find a countable family $P$ that satisfies
condition (i) of Lemma~\ref{lemma:suff} but not Lemma~\ref{lemma:necess}.

For each $k\in\N$ we define an infinite matrix
$\alpha^{(k)}=(\alpha^{(k)}_{ij})_{i,j\in\N}$
by setting
$$
\alpha^{(k)}_{ij}=\left\{
\begin{array}{ll}
ij, & i\le k\\
i, & i>k.
\end{array}
\right.
$$
Fix a bijection $\varphi\colon\N^2\to\N$ such that $\varphi(i,j)<\varphi(k,l)$
whenever $i+j<k+l$. For each $k\in\N$ define a sequence $p^{(k)}=(p^{(k)}_n)_{n\in\N}$
by $p^{(k)}_n=\alpha^{(k)}_{\varphi^{-1}(n)}$. Finally, set $P=\{ p^{(k)} : k\in\N\}$.

\begin{theorem}
The set $P$ has the following properties:
\begin{itemize}
\item[{\upshape (i)}] each sequence $p\in P$ has a bounded subsequence;
\item[{\upshape (ii)}] each sequence $\{ x_n\}\subset A(P)$ satisfying \eqref{bai}
is unbounded.
\end{itemize}
Therefore the algebra $A(P)$ has a locally bounded a.i., but does not have
a b.a.i.
\end{theorem}

Together with Theorem~\ref{thm:strflat2} this gives the following.

\begin{corollary}
There exists a commutative Fr\'echet-Arens-Michael algebra $A$
such that the trivial Fr\'echet $A$-module $\CC=A_+/A$ is strictly flat,
but $A$ does not have a b.a.i.
\end{corollary}

\section{Remarks on flat cyclic Fr\'echet modules\\
and amenable Fr\'echet-Arens-Michael algebras}
\label{sec:amen}

The second part of the Helemskii-Sheinberg theorem can be generalized as follows.
\begin{theorem}
\label{thm:flat}
Let $A$ be a Fr\'echet-Arens-Michael algebra, let $A=\varprojlim A_n$ be the
Arens-Michael decomposition of $A$, and let $I\subset A_+$ be a closed left
ideal. For each $n\in\N$, denote by $I_n\subset (A_n)_+$ the closure of the image
of $I$ under the canonical map $A_+\to (A_n)_+$. Suppose that $I_n$ is weakly
complemented in $(A_n)_+$ for each $n\in\N$.
Then the following conditions are equivalent:
\begin{itemize}
\item[{\upshape (i)}]
$A_+/I$ is flat;
\item[{\upshape (ii)}]
$A_+/I$ is strictly flat;
\item[{\upshape (iii)}]
$I$ has a right locally bounded a.i.
\end{itemize}
\end{theorem}

The condition {\em ``$I_n$ is weakly
complemented in $(A_n)_+$ for each $n\in\N$''} looks rather unnatural, but
we do not know how to put it into a more reasonable form. In particular,
we do not know the answers to the following questions:

\medskip\noindent
\textbf{Open problems.}
(1) Does the above condition depend on the choice of a sequence
of submultiplicative seminorms that gives the Arens-Michael decomposition of $A$?

(2) Suppose that $I$ is weakly complemented in $A_+$.
Does it follow that $I_n$ is weakly
complemented in $(A_n)_+$ for each $n\in\N$?

(3) Conversely, suppose that $I_n$ is weakly
complemented in $(A_n)_+$ for each $n\in\N$. Does it follow that
$I$ is weakly complemented in $A_+$?

(4), (5), (6). The same as (1), (2), (3), with ``weakly complemented''
replaced by ``complemented''.

\medskip
Fortunately, there is an important situation where the above-mentioned
difficulties disappear. First let us recall some standard notation \cite{X1}.
Given a Fr\'echet algebra $A$, the algebra opposite to $A$ is denoted by
$A^{\mathrm{op}}$. Set $A^e=A_+\Ptens A_+^{\mathrm{op}}$ and denote by
$\pi\colon A^e\to A_+$ the linear continuous map uniquely determined
by $a\otimes b\mapsto ab$. The kernel of this map is a complemented
closed left ideal of $A^e$. It is denoted by $I^\Delta$
and is called the {\em diagonal ideal} of $A^e$.

Recall that a Fr\'echet algebra $A$ is said to be {\em amenable} if
$A_+$ is a flat Fr\'echet $A$-bimodule (see \cite{X1}). Recall also that the category
of Fr\'echet $A$-bimodules is isomorphic to the category of left
unital Fr\'echet $A^e$-modules. The canonical morphism $\pi\colon A^e\to A_+$
determines an isomorphism between $A^e/I^\Delta$ and $A_+$ in $A^e\lmod$.
Therefore the question of whether or not $A$ is amenable is equivalent
to the question of whether or not the left cyclic Fr\'echet $A^e$-module
$A^e/I^\Delta$ is flat.

Now suppose that $A$ is a Fr\'echet-Arens-Michael algebra and $A=\varprojlim A_n$
is its Arens-Michael decomposition. Then we have $A^e=\varprojlim A^e_n$
(see, e.g., \cite{Kothe2}), and it is easy to see that the closure of the
canonical image of $I^\Delta$ in $A_n^e$ is precisely the diagonal ideal of $A^e_n$.
Since the latter is complemented in $A^e_n$, Theorem~\ref{thm:flat}
implies the following.

\begin{theorem}
Let $A$ be a Fr\'echet-Arens-Michael algebra. The following conditions
are equivalent:
\begin{itemize}
\item[{\upshape (i)}] $A$ is amenable;
\item[{\upshape (ii)}] $I^\Delta$ has a right locally bounded a.i.
\end{itemize}

If, in addition, $A$ is quasinormable,
then {\upshape (i)} and {\upshape (ii)} are equivalent to

\begin{itemize}
\item[{\upshape (iii)}] $I^\Delta$ has a right b.a.i.
\end{itemize}
\end{theorem}

\noindent\textbf{Open problems.}
(7) Does there exist an amenable Fr\'echet-Arens-Michael
algebra such that $I^\Delta$ does not have a right b.a.i.?

(8) Does there exist a quasinormable amenable Fr\'echet-Arens-Michael
algebra with this property?

\bibliographystyle{amsalpha}

\end{document}